
\input amssym.def

\magnification=1200
\baselineskip=14pt

\def\qed{{$\vrule height4pt depth0pt width4pt$}}
\def\ss{{\smallskip}}
\def\ms{{\medskip}}
\def\bs{{\bigskip}}
\def\ni{{\noindent}}

\def\s{{\sigma}}

\def\G{{\Gamma}}
\def\A{{\cal{A}}}

\def\Z{{\Bbb Z}}

\def\<{{\langle}}
\def\>{{\rangle}}

\def\m1{{\quad ({\rm mod\ 1})}}

\def\={\buildrel \cdot \over =}

\def\stab{{\rm Stab}}



\centerline{\bf An Invariant of Finite Group Actions on Shifts of Finite Type}
\smallskip

\centerline{Daniel S. Silver and Susan G. Williams} \ss

\ms
\footnote{} {The authors were partially supported by NSF grants
DMS-0071004 and DMS-0304971. The first author was partly supported by a grant from EGIDE, and the second
by a grant from CNRS.}
\footnote{}{2000 {\it Mathematics Subject Classification.}  
Primary 37B10; secondary 20F38, 57M27, 57R56.}

\noindent {\narrower\narrower\smallskip\noindent  {\bf
Abstract:}  We describe a pair of invariants for actions of finite groups on shifts of finite type, the left-reduced and right-reduced shifts.  The left-reduced shift was first constructed by U. Fiebig, who showed that its zeta function is an invariant, and in fact equal to the zeta function of the quotient dynamical system.  We also give conditions for expansivity of the quotient, and applications to combinatorial group theory, knot theory and topological quantum field theory.

\smallskip}
\ms


\ni {\bf 1. Introduction.} Let $(X,\s)$ be shift of finite type (see section 2 for definitions).   By an {\sl action} of a finite group $G$ on $X$ we mean a homomorphism $\phi$ of $G$ into the automorphism group of
$X$.    We say that $G$-actions $(X,\phi)$ and
$(Y,\psi)$  are  {\sl conjugate} if there is a topological conjugacy $\eta:X\to Y$ with
$\eta\circ\phi(g)=\psi(g)\circ \eta$ for all $g$ in $G$.  

Identifying points in the same $G$-orbit gives a quotient dynamical system $X/\phi$, with the quotient topology and homeomorphism induced by $\s$.  This
 need not be a shift of finite type, or even an expansive dynamical system.  In fact we show
in section 4 that for irreducible $X$,
$X/\phi$ is a shift of finite type if the quotient map is constant-to-one, and nonexpansive otherwise.  Ulf-Rainer
Fiebig {\bf[Fi93]} showed how to construct a shift of finite type that is an equal-entropy
factor of
$X$ and has the same number of period $n$ points, for every $n$, as the quotient
dynamical system $X/\phi$. In section 3 we examine Fiebig's shift and a mirror variant, which we call the left-reduced and right-reduced shifts $_\phi X$ and $X_\phi$ of the pair $(X,\phi)$.  Fiebig's construction depends, {\it a priori,} on the choice of a
presentation of
$X$ on which the action of $G$ is by one-block automorphisms.  A satisfying intrinsic
definition of the reduced shifts remains elusive.  However, we show in Theorem 3.6 that the conjugacy classes
of $_\phi X$ and $X_\phi$ are independent of the choice of presentation for $X$ and are thus  invariants of
the group action $(X,\phi)$.  Theorem 3.8 and Corollary 3.9 describe the behavior of these invariants under resolving and closing factor maps.

Our work was motivated by questions arising in applications of symbolic dynamics to
combinatorial group theory and topology.  We describe some of these applications in section 5.  In particular, a theorem of Patrick Gilmer [{\bf Gi99}] concerning a class of topological quantum field theories defined by Frank Quinn [{\bf Qu95}] is recovered as a special case of Theorem 3.6.

We thank Mike Boyle, Fran\c cois Blanchard and especially Bruce Kitchens for valuable
discussions.  We thank the referee for helpful comments.  We are grateful to the University of Maryland, the Institut de
Math\'ematiques de Luminy and the Centre de Math\'ematiques et Informatique, Marseille, for
their hospitality while much of this work was being done.\bs

\ni{\bf 2. Background.} In this section we briefly review some basic notions of symbolic dynamics.  For more background, includings zeta functions, Bowen-Franks groups and the technique of state-splitting used in the proof of Theorems 3.6 and 3.8, we refer the reader to  [{\bf LM95}] or [{\bf Ki98}].

Let $\A$ be a finite set or {\sl alphabet} of symbols.  The {\sl full shift} on $\A$ is the dynamical system consisting of the space $\A^\Z$ with the product topology and the left shift homeomorphism $\s$ given by $(\s  x)_i=(x_{i+1})$.  A {\sl subshift}  $X$ is a closed $\s$-invariant set of some full shift.  For $n\in {\Bbb N}$ the {\sl $n$-blocks} of $X$ are the symbol sequences of length $n$ that appear as consecutive entries of some $x\in X$.

Given a finite directed graph $\G =(V,E)$ with adjacency matrix $A$, the associated {\sl shift of finite type} (SFT)  $X_A$ is the subshift of $E^\Z$ consisting of bi-infinite sequences $(x_i)$ of edges that correspond to paths in $\G$. Vertices of $\G$, which form the index set of $A$, are also called {\sl states} of $X_A$. If there is an edge from state $i$ to state $j$ we say $j$ is a {\sl follower} of $i$, and $i$ a {\sl predecessor} of $j$. We will always assume that every state has a follower and predecessor; otherwise we can remove that state and its adjacent edges without changing $X_A$.  If $\G$ has no parallel edges, so that $A$ is a zero-one matrix,  the edge sequence  $(x_i)$ is determined by the sequence $(v_i)$ of initial vertices of these edges.  
The SFT $X_A$ is {\sl irreducible} if the graph $\G$ is strongly connected, that is, there is a path from any state to any other state.  (This is equivalent to topological transitivity of the dynamical system.) 

A  {\sl homomorphism} between dynamical systems $(X,S)$ and $(Y,T)$ is a continuous, map $\theta:X\to Y$ with $\, \theta\circ S=T\circ\theta$.  Epimorphisms are also called {\sl factor maps} and isomorphisms are {\sl (topological) conjugacies}.  If $\Theta$ is a map from the set of $n$-blocks of a shift space $(X,\s)$ to the alphabet of a shift space $(Y,\s)$ and $m\in\Z$, we can define a homomorphism $\theta: X\to Y$ by $\theta((x_i))=(y_i)$ where $y_i=\Theta(x_{i-m},\ldots,x_{i-m+n-1})$; $\eta$ is an {\sl $n$-block map} with memory $m$.  For convenience we sometimes use the same symbol for the map on $n$-blocks and the map on $X$; when we define a map via its action on blocks we will take $m=0$.  Every homomorphism between subshifts is an $n$-block map for some $n$.

A one-block map $\theta: X\to Y$ is {\sl right-resolving} if whenever $ab$ and $ac$ are 2-blocks of $X$ with $\Theta(b)=\Theta(c)$ we have $b=c$.  If $X$ is a shift of finite type described by a graph $\G$ with no parallel edges, this is equivalent to the condition that edges with the same initial vertex have distinct images under $\Theta$.  A homomorphism of subshifts is {\sl right-closing} if it does not identify left-asymptotic points, i.e. points $(x_i)$ and $(x'_i)$ with $x_i=x'_i$ for $i\leq N$. Right-resolving maps are right-closing.  Right-closing factor maps are bounded-to-one and hence preserve topological entropy.  Left-resolving and left-closing maps are defined analogously.  One-block maps that are both left and right-resolving are {\sl bi-resolving} and homomorphisms that are left and right-closing are {\sl bi-closing}.  A factor map between irreducible shifts of finite type is constant-to-one if and only if it is bi-closing [{\bf Na83}].

\bs

\ni{\bf 3. The reduced shift of a group action.} In what follows we will consider $G$-actions $(X,\phi)$, always assuming that $X$ is a
shift of finite type and $G$ is a finite group. When there is no danger of confustion, we will  denote the
image of a point $x$ under $\phi(g)$ by $gx$, and its orbit under $G$ by $Gx$.
Consider the special case of a group action $(X_A,\phi)$ where $A$ is an $n\times n$ matrix over $\{0,1\}$, $X_A$ is the associated shift of finite type, and each
$\phi(g)$ is a one-block automorphism induced by a permutation of the states of $X_A$.  We will
call such an action a {\sl permutation action}.  We write $gi$ for the image of state $i$ under the permutation $\phi(g)$, and $Gi$ for its
orbit.\bs

\ni{\bf Proposition 3.1.} [{\bf AKM85}] Every finite group action on a SFT is conjugate to a
permutation action.\bs

\ni{\bf Definition 3.2.} Given a permutation action $(X_A,\phi)$, we define the associated {\sl
right-reduced shift} to be the SFT $X_\phi$ given by the matrix $A_\phi$ such that its states are
$G$-orbits of states of $X_A$ and $A_\phi(Gi,Gj)=\sum_{k\in Gj}A(i,k).$  This is well defined
since if
$i'=gi$ then $$\sum_{k\in Gj}A(i',k)=\sum_{k\in Gj}A(gi,gk)=\sum_{k\in Gj}A(i,k).$$
Analogously, the {\sl left-reduced shift} $_\phi X$ is given by the matrix $_\phi A$ with the same state space and $_\phi A(Gi,Gj)=\sum_{k\in Gi}A(k,j).$\bs

The  left-reduced shift was constructed (but not named) by U.-R. Fiebig [{\bf Fi93}], who proved the rather surprising result that it has  the same zeta function as $X_A/\phi$, or equivalently, the same number of period $n$ points for every $n$. A
permutation action on $X_A$ is also a permutation action on the inverse shift given by the transpose $A^t$, and it is easy to see that
the left-reduced shift for $(X_{A^t},\phi)$ is the inverse of the right-reduced shift for
$(X_A,\phi)$. Thus the left-reduced and right-reduced shifts have analogous properties, although
Example 3.7 shows that they need not be either conjugate or inverse conjugate.  For simplicity, we state and prove most of our results for the right-reduced shift, which we will refer to as the {\sl reduced shift} when there is no danger of confusion.  
(This dextrocentric preference 
follows a tradition in the literature of shifts of finite type, and is also in convenient agreement with the work of Quinn and Gilmer that we discuss in Section 5.)  

The construction of the reduced shift, despite its simplicity, is not natural in the sense that it relies on the presentation of the group action as a permutation action.
We can easily describe a  right-resolving factor map $\eta: X_A \to X_\phi$:  for each pair of states $i$, $j$ of $X_A$, pick a bijection from the set of edges $(i,k)$ with $k\in Gj$ to the set of edges in $X_\phi$ from $Gi$ to $Gj$.  However, different choices may give nonisomorphic factor maps.  (Factor maps $\eta$, $\theta$ from $X$ to $Y$ are isomorphic if there are automorphisms $\alpha$, $\beta$ of $X$ and $Y$ respectively with $\theta\circ\alpha=\beta\circ\eta$.) In general, neither $X_\phi$ nor $X/\phi$ appears in a natural way as a factor of the other.

\bs
\ni {\bf Example 3.3.} 
Let $X$ be the full shift on the elements of the symmetric group $S_3$.  Let $G=S_3$ act on the states of $X$ by conjugation, so that $g\in G$ takes the state $v$ to $g^{-1}vg$.  This action induces  a permutation action $\phi$ of $G$ on $X$.  The states of the reduced shift $X_\phi$ are the three conjugacy classes of elements of $S_3$.  If we list these in the order $[(1)], [(12)], [(123)]$ then $X_\phi$ is given by the matrix
$$A_\phi=\pmatrix{1&3&2\cr 1&3&2\cr1&3&2\cr}.$$ 

There are many ways to define a 2-block right-resolving factor map from $X$ to $X_\phi$.  No choice respects the action of $G$.  For example, suppose the two-block $(12), (23)$ is sent to an edge $e$ of $X_\phi$ ($e$ must be one of the self-loops on the vertex [(12)].)  To respect the conjugation by (13) we would need to send the 2-block $(23),(12)$ to $e$ as well, while to respect the conjugation by (123) we would need to send $(23),(31)$ to $e$.  It is easy to see that different choices can give non-isomrophic factor maps: for example, we may send the fixed points of $X$ to distinct fixed points of $X_\phi$, or make some identifications.
\bs
\ni {\bf Example 3.4.} Let $$A=\pmatrix{1&1&1\cr 1&1&0\cr 1&0&1}$$
and let $\phi$ be as in the preceding example. Then $X_\phi$ is given by the graph with symbolic matrix
$$\tilde A_\phi=\pmatrix{a&b+c\cr d&e}.$$  Let $q$ denote the quotient map from $X$ to $X/\phi$.   There is no right-resolving factor map $\eta:X\to X_\phi$ satisfying $\eta=\xi\circ q$ for some homeomorphism $\xi:X/\phi\to X_\phi$, since $\eta$ could not identify the left-asymptotic points $1^\infty2^\infty=(\ldots,1,1,2,2,\ldots)$ and $1^\infty3^\infty$ which are identified by $q$.  There is also
 no factor map $\eta:X\to X_\phi$ satisfying $q=\xi\circ \eta$ for some homomorphism $\xi:X_\phi\to X/\phi$.  For $\eta$ would have to take fixed points $1^\infty$, $2^\infty$, $3^\infty$ to $a^\infty$, $e^\infty$, $e^\infty$ respectively.  Since $\eta$ is an $n$-block map for some $n$, it would have to identify  $2^\infty1^n 2^\infty$ and $2^\infty1^n 3^\infty$, which are not identified by $q$.\bs

We may obtain $A_\phi$ from $A$ as a matrix product.  Suppose $A$ has $n$ states and $\bar n$
$G$-orbits of states.  Let $V_\phi$ be the $n\times \bar n$ matrix with
$V_\phi(i,Gi)=1$ for all $i$ and the remaining entries 0.  Take $U_\phi$ to be any left
inverse of $V_\phi$. Thus $U_\phi$ is $\bar n\times n$ and $U_\phi(Gi,i_0)=1$ for a selected
representative $i_0$ of the orbit $Gi$.  One verifies easily that $A_\phi=U_\phi A V_\phi$.  Letting
$P_{\phi(g)}$ denote the matrix of the permutation $\phi(g)$, we have
$P_{\phi(g)}V_\phi =V_\phi$ for all $g\in G$.

An {\sl elementary strong shift equivalence} between square non-negative integral matrices $A$
and $B$ consists of a pair $(R,S)$ of non-negative integral matrices (not necessarily square) with $RS=A$ and
$SR=B$.  Matrices $A$ and $B$ are {\sl strong shift equivalent} (SSE) if they are linked by a
chain of elementary strong shift equivalences, that is, if there are  square nonnegative integral matrices $A=A_0, A_1,\ldots A_n=B$ with an elementary strong shift equivalence  from $A_{i-1}$ to $A_i$ for $1=1,\ldots ,n$. The Decomposition Theorem of R.F. Williams
[{\bf Wi73}] states that SFT $X_A$ and $X_B$ are conjugate if and only if $A$ and $B$ are SSE.

To an elementary SSE $(R,S)$ of 0-1 matrices $A$, $B$ we may canonically associate a
conjugacy from $X_A$ to $X_B$, as described in [{\bf Ki98}], Lemma 2.1.16, or [{\bf LM95}],
p.228.  Briefly, since $A=RS$, for each nonzero entry $(i,j)$ of $A$ there is a unique state $k$ of $B$ with $R(i,k)=S(k,j)=1$.  This allows us to associate to each $x\in X_A$ a bi-infinite sequence of states of $B$.  The identity $SR=B$ implies that there is a unique edge sequence $y\in X_B$ connecting this sequence of states.  The map $\eta(x)=y$ is a conjugacy.
We observe that this conjugacy will carry a permutation action $\phi$ of $G$ on $X_A$ to a
permutation action $\psi$ of $G$ on $X_B$ if and only if 
$RP_{\psi(g)}=P_{\phi(g)}R$ and $SP_{\phi(g)}=P_{\psi(g)}S$ for all $g\in G$.  For notational simplicity we will view $X_A$ and $X_B$ as having disjoint state spaces $I$ and $I'$ that each admit an action of $G$ by permutations, so that $gi$ denotes $\phi(g) i$ for $i\in I$ and $\psi(g) i$ for $i\in I'$.\bs

\ni{\bf Lemma 3.5.} Let $(R,S)$ be an elementary SSE of 0-1 matrices $A$ and $B$ that induces a
conjugacy of permutation actions $\phi$, $\psi$ on $X_A$, $X_B$ respectively.  Then the
reduced shifts $X_\phi$ and $X_\psi$ are conjugate. \bs 

\ni{\bf Proof.} We claim that the pair $(U_\phi RV_\psi, U_\psi SV_\phi)$ is an elementary SSE
between $A_\phi$ and $B_\psi$. For each $g\in G$,
$$P_{\psi(g)}SV_\phi =SP_{\phi(g)}V_\phi =SV_\phi,$$
which implies that the $i$-th and $gi$-th rows of $SV_\phi$ are identical for all $i$.  Now,
$V_\psi U_\psi$ has 1 in the $(i,i_0)$ entry where $i_0$ is the distinguished representative of
the orbit $Gi$, and 0 elsewhere.  Hence $V_\psi U_\psi SV_\phi = SV_\phi$.  This yields
$$(U_\phi RV_\psi)(U_\psi SV_\phi) = U_\phi RSV_\phi = U_\phi AV_\phi = A_\phi.$$
Similarly, $(U_\psi SV_\phi)(U_\phi RV_\psi)=B_\psi$. \qed\bs

\ni{\bf Theorem 3.6.} If $(X_A, \phi)$ and $(X_B,\psi)$ are conjugate permutation actions on
SFT then the reduced shifts $X_\phi$ and $X_\psi$ are conjugate.\bs

\ni{\bf Proof.} By the Decomposition Theorem, the conjugacy can be expressed as a composition
of conjugacies corresponding to elementary SSE.  In light of the preceding lemma, it suffices
to show that the decomposition can be carried out in such a way that each elementary SSE
induces a conjugacy of permutation actions.

We follow the proof of the Decomposition Theorem in [{\bf LM95}] (Theorem 7.1.2). By
passing from
$X_A$ by a chain of elementary SSE to a higher block presentation we may replace the
original conjugacy by a 1-block map.  The  $n$-block presentation of $X_A$ is again given by a 0-1 matrix, and each permutation $\phi(g)$ of the states of $X_A$ naturally induces a permutation
of $n$-blocks that is a conjugate permutation action on the $n$-block presentation.

Thus we may assume we have a 1-block conjugacy $\eta:X_A\to X_B$.  If $\eta^{-1}$ is also a
1-block map then the conjugacy is simply a renaming of symbols, and the result is clear.  If
$\eta^{-1}$ has anticipation $k>1$ we use an out-splitting of the graph of $X_A$ to reduce the
anticipation.    At each vertex of the graph of $X_A$ we partition the outgoing edges according to
their images under $\eta$.  These partition elements are the states of a conjugate shift
$X_{\tilde A}$, where ${\tilde A}$ is a 0-1 matrix, and $\eta$ induces a one-block conjugacy
$\tilde\eta$ from $X_{\tilde A}$ to the two-block presentation of  $X_B$ such that $\tilde\eta^{-1}$ has the same memory as
$\eta$ but anticipation $k-1$.  Since
$\eta$ intertwines the $G$-actions on $X_A$ and $X_B$, if two edges leaving state $i$
have the same image under $\eta$ then any $g\in G$ carries them to edges leaving state
$gi$ with the same image under $\eta$.  This determines a permutation action on $X_{\tilde A}$.
It is easy to see that all of these conjugacies preserve the $G$-actions.

The memory of $\eta^{-1}$ may be reduced by in-splittings in an analogous fashion.  An
induction argument finishes the proof. \qed\bs

In view of Proposition 3.1 and Theorem 3.6, we can speak of the (left or right) reduced shift of an arbitrary $G$-action on a SFT with the understanding that it is well defined up to topological conjugacy.

We say the $G$-action $(Y,\psi)$ is a {\sl factor} of the $G$-action $(X,\phi)$ if there is a
factor map $\eta:X\to Y$ with $\eta\circ \phi(g)=\psi(g)\circ \eta$ for all
$g\in G$.  In this case the quotient dynamical system $Y/\psi$ is a factor of $X/\phi$, and it is natural to ask whether the reduced shift $Y_\psi$ is a factor of $X_\phi$.  Example 3.7 shows that this can fail even for almost invertible factor maps between irreducible SFT.  Theorem 3.8 and its corollary give a positive answer for right-resolving and right-closing factor maps.  Analogous results hold for left-closing and left-resolving factors.
\bs

\ni {\bf Example 3.7.} Let $$A=\pmatrix{1&0&1&0&1&0\cr 0&1&0&1&0&1\cr 1&1&1&0&0&0\cr
1&1&0&1&0&0\cr 1&1&0&0&1&0\cr 1&1&0&0&0&1\cr}$$
and let $\phi$ be the permutation action on $X_A$ of the cyclic group $G\cong{\Z}/4$
generated by the permutation $(12)(3456)$. Then 
$$A_\phi=\pmatrix{1&2\cr 2&1},\ \ _\phi A=\pmatrix{1&1\cr 4&1}.$$
Here we have taken the states of the reduced shifts to be $G1=\{1,2\}$ and $G3=\{3,4,5,6\}$ in
that order.  The left-reduced shift is conjugate to neither the right-reduced shift nor its
inverse, as they have non-isomorphic Bowen-Franks groups ${\Z}^2/(I-A_\phi){\Z}^2 \cong{\Z}/2 \oplus {\Z}/2$ and ${\Z}^2/(I- \!_\phi A){\Z}^2 \cong {\Z}/4$.

Let $$B=\pmatrix{ 1&1&1&1&1\cr 1&1&0&0&0\cr
1&0&1&0&0\cr 1&0&0&1&0\cr 1&0&0&0&1\cr}.$$
A one-block factor map of $X_A$ onto $X_B$ is obtained by identifying the first two
states of $A$, and this induces a $\Z/4$-action $\psi$ on $X_B$ that cyclically permutes the last  four states of $B$.  Points of $X_B$ that are right-asymptotic to the fixed point on the first state have two preimages in $X_A$ but every other point has a unique preimage.  We have 
$$B_\psi=\pmatrix{1&4\cr 1&1}.$$
The right-reduced shift $X_\psi$ is not a factor of $X_\phi$, as can be seen from either Theorem
4.2.16 or Theorem 4.2.19 of [{\bf Ki98}].  (The latter theorem says that a factor map between
irreducible SFT induces an epimorphism of their Bowen-Franks groups.)

A simpler, but reducible, example of these phenomena is given by
$$A=\pmatrix{1&0&1\cr 0&1&1\cr 0&0&1\cr}$$
with permutation action of $G\cong{\Z}/2$
generated by the permutation $(12)$.   Then 
$$A_\phi=\pmatrix{1&1\cr 0&1},\ \ _\phi A=\pmatrix{1&2\cr 0&1}.$$
Let $B$ be the quotient $_\phi A$ with trivial $G$ action $\psi$.  Then $A_\phi$ and $_\phi A$ are nonconjugate, and $B_\psi=B$ is not a quotient of $A_\phi$.\bs

\ni {\bf Theorem 3.8.} Suppose the permutation action $(X_B,\psi)$ is a factor of the permutation action $(X_A,\phi)$ by a right-resolving 1-block map $\eta$.  Then there are right resolving 1-block maps $\bar\eta:X_\phi\to X_\psi$, $\theta_1:X_A\to X_\phi$ and $\theta_2:X_B\to X_\psi$ with $\theta_2\circ\eta=\bar\eta\circ\theta_1$.
\bs

\ni{\bf Proof.} The 1-block map $\eta$ induces a graph homomorphism from the graph of $X_A$ to the graph of $X_B$, which takes $i$ to a vertex we denote by $\eta(i)$.  Since $\eta$ is right-resolving, distinct edges with the same initial vertex $i$ have distinct images under this graph homomorphism.  Since $\eta$ preserves the group action, the images of vertices in an orbit $Gj$  comprise the orbit $G\eta(j)$.  Hence $\eta$ gives a bijection from the set of edges with initial vertex $i$ and terminal vertex in $Gj$ to the set of edges with initial vertex $\eta(i)$ and terminal vertex in $G\eta(j)$.  This means the number of edges in the graph of $X_\phi$ from $Gi$ to $Gj$ is equal to the number of edges in the graph of $X_\psi$ from $G\eta(i)$ to $G\eta(j)$; we define the 1-block map $\bar\eta$ by making any choice of bijections.  We may also let $\theta_2$ be given by any choice of bijections from the edges in the graph of $X_B$ with initial state $i'$ and terminal state in $Gj'$ to the edges in the graph of $X_\psi$ from state $Gi'$ to state $Gj'$. 
There is now a unique 1-block map $\theta_1$ which gives the desired commutativity.

Because $\eta$ and $\theta_2$ are right-resolving, the composite map $\theta_2\circ\eta=\bar\eta\circ\theta_1$ is also right-resolving.  Now since $\theta_1$ is onto, it is easy to see that $\bar\eta$ must be right-resolving as well. \qed\bs

\ni{\bf Corollary 3.9.} If the $G$-action $(Y,\psi)$ is a factor of the $G$-action $(X,\phi)$ by a right-closing map $\eta$ then $Y_\psi$ is a right-closing factor of $X_\phi$.
\bs

\ni{\bf Proof.}  By Proposition 3.1, we can assume from the start that the actions are permutation actions.  By Proposition 1 of [{\bf BKM85}], there is a topological conjugacy $\pi: \tilde X \to X$ such that $\eta\circ\pi$ is a right-resolving 1-block map.   It is easy to see from their construction that $\pi^{-1}$ induces a permutation action on $\tilde X$. We apply Theorem 3.8 to $\eta\circ\pi$, then use Theorem 3.6 together with the fact that the composition of a right-resolving map with a topological conjugacy is right-closing. \qed\bs


\ni{\bf 4. $G$-stabilizers and $G$-orbits of periodic points.}  
Given a $G$-action $(X, \phi)$, for $x\in X$ we denote by $\stab (x)$ the $G$-stabilizer of $x$, that is, the subgroup of all $g\in G$ with $gx=x$.  For permutation actions, the $G$-stabilizer of a state or word of $X$ may be defined similarly.  In this case the $G$-stabilizer of $x=(x_i)$ is $\stab (x) = \bigcap _{i\in \Z} \stab( x_i)$.
The number of preimages of a point $Gx\in X/\phi$ under the quotient map is the index in $G$ of $\stab(x)$.\bs

\ni{\bf Theorem 4.1.} Let $(X,\phi)$ be a $G$-action on an irreducible shift of finite type. (i) If every $x\in X$ has the same stabilizer then the quotient map is constant-to-one, and $X/\phi$, $X_\phi$ and $_\phi X$ are all conjugate shifts of finite type.  (ii) If some pair of points of $X$ have different $G$-stabilizers then $X/\phi $ is nonexpansive.\bs

\ni{\bf Proof.} We may assume $(X,\phi)$ is a permutation action.  Suppose first that every point has the same $G$-stabilizer $H$. Since every state of $X$ appears in some periodic point of period at most $n$, the number of states of $X$, by passing to the $n$-block presentation of $X$ we can assume every state of $X$ has stabilizer $H$.  Now it is clear that the one-block map $i \mapsto Gi$ is a bi-resolving map from $X$ to $_\phi X=X_\phi$, and the image is topologically conjugate to $X/\phi$. 

Now suppose there are points of $X$ with different stabilizers. We first show there must be periodic points $u^\infty$, $v^\infty$ with different stabilizers.  Since $X$ is irreducible there is a periodic point $v^\infty$ containing all states, so that $\stab(v^\infty)$ is a subgroup of $\stab(x)$ for all $x\in X$.  Choose any $y\in X$ with $\stab(y)\neq \stab(v^\infty)$. Then $y$ contains a word $u$ such that $u^\infty\in X$, and $\stab(y)$ is a subgroup of $\stab(u^\infty)$.

Let $g\in \stab(u) \setminus\stab(v)$.  We can find words $w, w'$ such that $uw v$ and $vw'u$ are words of $X$.   Then (gu)(gw)(gv)=u(gw)(gv) is also a word of $X$.
For each positive integer $m$ set
$$\eqalign{ x^{(m)} &=v^{\infty}w'.u^{2m+1}wv^\infty \cr
y^{(m)} &=v^{\infty}w'.u^{2m+1}(gw)(gv)^\infty .\cr}
$$
(Here the point precedes the 0-coordinate.)
Then $x^{(m)}$ and $y^{(m)}$ have different images $[x^{(m)}]$, $[y^{(m)}]$ under the quotient map.  However, for every $n\in \Z$ the central $(2m+1)$-block of $\s^ny^{(m)}$ agrees with the central $(2m+1)$-block of either $\s^nx^{(m)}$ or $\s^ngx^{(m)}$, so that $[\s^nx^{(m)}]$ and $[\s^ny^{(m)}]$ are close in the quotient topology. Hence there is no expansive constant for $X/\phi$. \qed\bs

We next give a formula for the number of $G$-orbits of period $n$ points of $X$.  Note that
this is different from the number of period $n$ points of the quotient dynamical system
$X/\phi$, since the quotient map may change the period of a point.  An application of this result appears as Propositon 5.1 below.

For each $g\in G$ the set Fix$(g)$ of points of $X$ fixed by $g$ is again a SFT.  If we assume
that $\phi$ is a permutation action on $X_A$ then Fix$(g)=X_{A_g}$ where $A_g$ is the
principal submatrix of $A$ corresponding to the set of symbols fixed by $g$. (If this set is empty we take $A_g=(0)$.) \bs

\ni{\bf Proposition 4.2.} Let  $G$ act by permutations on a shift of finite
type $X_A$.  The number $N_n$ of
$G$-orbits of period $n$ points of $X_A$ is given by
$$N_n={1\over |G|}\sum_{g\in G} {\rm trace} (A_g^n).$$
Hence the sequence $\{N_n\}$ satisfies a linear homogeneous recurrence relation with constant coefficients.\bs
 
\ni {\bf Proof.} A combinatorial result of Cauchy and Frobenius commonly known as
the Burnside Lemma (cf.~[{\bf DM96}], p.24) says that if a group
$G$ acts on a set $S$ then the number of $G$-orbits of $S$ is 
$${1\over |G|}\sum_{g\in G} |{\rm Fix}(g)|,$$
where ${\rm Fix}(g)$ is the set of points fixed by $g$.  If we take $S$ to be the set of
period $n$ points of $X_A$ then ${\rm Fix}(g)$ is the set of period $n$ points of $X_{A_g}$,
which has cardinality trace$(A_g^n)$.  The sequence $\{{\rm trace} (A_g^n)\}$ satisfies the
linear recurrence with characteristic polynomial $\det(I-tA_g)$, so the sum satisfies the
recurrence relation given by the least common multiple of these polynomials.
\qed\bs


\ni{\bf 5. Representation shifts, knots and topological quantum field theories.}  In this section we describe a class of SFT called representation shifts that admit a natural group action, and outline applications to knot theory and topological quantum field theory.

Assume
that $\Pi$ is a finitely presented group with epimorphism $\chi:
\Pi\to \Z$, and let $x \in \chi^{-1}(1)$. Then $\Pi$ can be described as an {\it HNN extension}
$\langle x, B \mid x^{-1}ax = \phi(a),\ \forall a \in U\rangle$.
Here  $B$ is a finitely generated
subgroup of $K = {\rm ker}\ \chi$, and the map $\phi$ is an
isomorphism between finitely generated subgroups $U,V$ of $B$.
The subgroup $B$ is an {\it HNN base},
$x$ is a {\it stable letter}, $\phi$ is an {\it amalgamating
map}. Details can be found in [{\bf LS77}].

Conjugation by $x$ induces an automorphism  of $K$. Letting
$B_j = x^{-j}Bx^j, U_j = x^{-j}Ux^j$ and $V_j= x^{-j}Vx^j, j
\in \Z$, we can express $K$ as an infinite amalgamated free
product
$$K = \langle B_j\mid V_j = U_{j+1},\ \forall j \in
\Z\rangle.$$

For any finite group $G$, the set ${\rm Hom}(K,G)$ of representations of $K$ in $G$ may be viewed as a SFT.  The state set is ${\rm Hom}(U,G)$; an edge is an element $\rho_0\in{\rm Hom}(B,G)$ with initial state $\rho_0|_U$ and terminal state $\rho_0|_V\circ\phi$.  (Note that the cardinality of the edge set is bounded by $|G|^m$ where $m$ is the cardinality of a generating set of $B$.)  If $\rho=(\rho_j)$ is a bi-infinite path then the representations from $B_j$ to $G$ given by $y\mapsto \rho_j(x^jyx^{-j})$ have a unique common extension to an element of ${\rm Hom}(K,G)$ that we will also denote by $\rho$.  

We call this SFT the {\it representation shift} of $K$ in $G$ and denote it by $\Phi_G=\Phi_G(\Pi,\chi, x)$ (see [{\bf SW96}]). The usual topology on the SFT coincides with the compact-open topology on ${\rm Hom}(K,G)$, and the shift map $G$ can be described by $\s\rho(y)=x^{-1}yx$.  It is clear from this intrinsic description that the topological conjugacy class of $\Phi_G$ is independent of the choice of HNN base $B$. 

The group $G$ acts on $\Phi_G$ by inner automorphism of the image space, $(g\rho)(x)=g^{-1}\rho(x)g$.  The case where $G$ is the symmetric group $S_n$ is of particular interest in studying finite-index subgroups of $K$ (see [{\bf SW96}], [{\bf SW99$'$}] and [{\bf SW04}]). In this case the state space of the representation shift typically grows very quickly with $n$.  Since the group of inner automorphisms is isomorphic to $S_n$, the corresponding reduced shift will be considerably simpler.  (For $n\neq 6$ all automorphisms of $S_n$ are inner: see theorem 6.20 of [{\bf Is94}]).

Much of our original motivation came from the study of knots. A {\it knot} $k$ is a smoothly embedded circle in the $3$-sphere $\Bbb S^3$. For convenience, we assume that $k$ is oriented. Two knots are regarded as the same if they are ambiently isotopic.  Although the {\it knot group} $\Pi= \pi_1(\Bbb S^3 \setminus k)$ is essentially a complete invariant, unlocking all of its information is not a reachable task at present. Fortunately, many tractable invariants can be computed from the group. 

Let $x \in \Pi$ be the element represented by a meridian curve encircling $k$ with linking number $1$. The abelianization of $\Pi$ is infinite cyclic, and we let  $\chi: \Pi \to \Bbb Z$ be the abelianization homomorphism that maps $x$ to $1$.  The kernel $K$ of $\chi$ is the commutator subgroup of $\Pi$. Given a finite group $G$, the representation shift  $\Phi_G=\Phi_G(\Pi,\chi, x)$ is an invariant of $k$ [{\bf SW96}]  (see also [{\bf SW99}]). 

For any positive integer $n$, the set ${\rm Fix}(\s^n)$ of period $n$ points coincides with the set 
${\rm Hom}(\pi_1 M_n, G)$, where $M_n$ is the $n$-fold cyclic cover of $\Bbb S^3$ branched over $k$. Details can be found in [{\bf SW99}].  Invariants of branched covers $M_n$ are invariants of $k$. (Such invariants were first considered by J. Alexander and G. Briggs in the early 1900's, and they remain an important class.) The group $G$ acts on ${\rm Hom}(\pi_1 M_n, G)$ by inner automorphism as above, and the orbit set ${\rm Hom}(\pi_1 M_n, G)/G$ can be identified with the set of flat $G$-bundles on $M_n$. Proposition 4.2 immediately yields the following. \bs

\ni {\bf Proposition 5.1.} For any knot $k$ and finite group $G$, the number of flat $G$-bundles over the $r$-fold cyclic cover $M_n$ of $\Bbb S^3$ satisfies a linear recurrence. \bs

Proposition 5.1 should be compared to Theorem 4.2 of [{\bf SW99}] or Proposition  2.3 of [{\bf Gi99}], either of which shows that $|{\rm Hom}(\pi_1 M_n, G)|$ satisfies a linear recurrence. 

Gilmer's paper [{\bf Gi99}] is concerned with  topological quantum field theories (TQFT), and was another source of motivation for us. TQFT, which arose from quantum physics, offer a framework for the understanding invariants such as the Jones polynomial for links or Donaldson invariants of $4$-manifolds or the discovery of new invariants. 

Roughly speaking, a (d+1)-dimensional TQFT assigns to a $d$-dimensional oriented manifold $Y$, called a {\it space}, a module $Z(Y)$ over a coefficient ring $R$. When spaces $Y_1$ and $Y_2$ are the ``incoming" and ``outgoing" boundaries of a {\it spacetime}, an oriented $(d+1)$-dimensional manifold $X$, a homomorphism $Z_X: Z(Y_1) \to Z(Y_2)$ is assigned. In particular, we require $Z_{Y\times [0,1]} = {\rm id}: Z(Y) \to Z(Y)$, which implies that if $Z_X$ is nontrivial, then $X$ is topologically nontrivial (that is, not a product of a $d$-dimensional manifold with the unit interval). 
Various other axioms are imposed. For example,  $Z(Y_1 \sqcup Y_2) = Z(Y_1) \otimes_R Z(Y_2)$, where $\sqcup$ denotes disjoint union. Also, if $X_1$ has incoming (resp. outgoing) boundaries $Y_1$ (resp. $Y_2$) while $X_2$ has incoming (resp. outgoing) boundaries $Y_2$ (resp. $Y_3$):
$$Y_1\ {\buildrel X_1 \over \longrightarrow}\ Y_2\ {\buildrel X_2 \over \longrightarrow}\ Y_3,$$
then the associated homomorphisms compose in a natural way:
$$Z_{X_1 \cup_{Y_2} X_2} = Z_{X_2} Z_{X_1}.$$
For a more complete discussion, the  reader might consult [{\bf Qu95}] or [{\bf At89}]. In Quinn's very general approach, manifolds can be replaced by finite CW complexes with variously defined boundaries. 

Quinn [{\bf Qu95}] uses a finite group $G$ to construct a TQFT. For the sake of simplicity, we describe his TQFT only for a special case that arises in knot theory.  As above, let $k$ be an oriented knot with group $\Pi = \pi_1(\Bbb S^3 \setminus \ell)$ and abelianization homomorphism $\chi: \Pi \to \Bbb Z$ mapping the distinguished element $x \in \Pi$ to $1$. 
It is well known that $ {\Bbb S}^3 \setminus k$ admits a smooth map $f$ to ${\Bbb S}^1$ inducing $\chi$ on first homology groups and such that the preimage of a regular value is a connected orientable surface,  the interior of a surface $Y \in {\Bbb S}^3$ with boundary $k$, called a {\it Seifert surface} for the knot. Cutting ${\Bbb S}^3$ along $Y$ produces a compact manifold $X$ with two boundary components $Y_1, Y_2$ that are copies of $Y$. 

Quinn's TQFT assigns a ${\Bbb Q}$-vector space $Z(Y_i)$ to $Y_i$, $i=1, 2$, with basis consisting 
of homomorphisms from $\pi_1(Y_i)$ to $G$ modulo inner automorphisms of $G$; in other words, $Z(Y_i)$ has basis ${\rm Hom}(\pi_1(Y_i), G)/G$. Choosing basepoints $y_i \in Y_i$ and a path $s$ in $X$ connecting $y_1$ and $y_2$, one defines a homomorphism $s_*: \pi_1(Y_2, y_2) \to \pi_1(X, y_1)$ sending a loop $\gamma$ at $y_2$ to $s\gamma s^{-1}$. For any $\beta: \pi_1(X, y_1)\to G$, let $\beta_1: \pi_1(Y_1, y_1) \to G$  be the composition of the map $\pi_1(Y_1, y_1) \to \pi_1(X, y_1)$ induced by inclusion, and $\beta$.  Let  $\beta_2: \pi_1(Y_2, y_2) \to G$ be the composition of $s_*$ and $\beta$.  Combining these ingredients,  we define $Z_{Y_1, Y_2}$ to be the homomorphism from $Z(Y_1)$ to $Z(Y_2)$ that maps $[\alpha] \in {\rm Hom}(\pi_1(Y_1), G)/G$ to  $\sum_\beta [\beta_2],$
where the sum is taken over all $\beta: \pi_1(X, y_1) \to G$ such that $ \beta_1 = \alpha.$
In a sense, $Z_{Y_1, Y_2}$ records the various ways that $\beta_2$ extends over
$\pi_1(X)$ modulo inner automorphisms of $G$. 

We identify $Z_{Y_1, Y_2}$ with its matrix representation (with respect to the given bases). Since its entries are nonnegative integers, it defines a shift of finite type. In fact, it is the reduced shift of the representation shift $\Phi_G=\Phi_G(\Pi,\chi, x)$, defined above. Since the representation shift is well defined, it follows from Theorem 3.6 that the strong shift equivalence class of $Z_{Y_1, Y_2}$ is an invariant of the knot $k$. This fact was established by Gilmer in [{\bf Gi99}] using topological methods.


\bs

\centerline{{\bf References}}
\bs
\item{[{\bf AKM85}]} R.L. Adler, B. Kitchens and B.H. Marcus, Finite group actions on shifts
of finite type, {\sl Ergod.\ Th.\ \& Dynam.\ Sys.\ \bf 5} (1985), 1--25.\ss

\item{[{\bf At89}]} M.F. Atiyah, Topological Quantum Field Theories, {\sl Publ.\ Math.\ Inst.\ Hautes\ Etudes\ Sci.\ \bf 68} (1989), 175--186. \ss

\item{[{\bf BKM85}]} M. Boyle, B. Kitchens and B. Marcus, A note on minimal covers for sofic systems, {\sl Proc.\ Amer.\ Math.\ Soc.\ \bf 95} (1985), 403--411.\ss

\item{[{\bf DM96}]} J.D. Dixon and B. Mortimer, {\sl Permutation Groups,} Springer-Verlag,
1996.\ss

\item{[{\bf Fi93}]} U. Fiebig, Periodic points and finite group actions on shifts of finite
type, {\sl Ergod.\ Th.\ \& Dynam.\ Sys.\ \bf 13} (1993), 485--514.\ss

\item{[{\bf Gi99}]} P.M. Gilmer, Topological quantum field theory and strong shift
equivalence, {\sl Canad.\ Math.\ Bull. \bf 42} (1999), 190--197.\ss

\item{[{\bf Is94}]} I.M. Isaacs, {\sl Algebra: A Graduate Course,} Brooks/Cole, Belmont, CA, 1994. \ss

\item{[{\bf Ki98}]} B.P. Kitchens, {\sl Symbolic Dynamics: One-sided, two-sided and countable
state Markov chains,} Springer-Verlag, Berlin 1998.\ss

\item{[{\bf LM95}]} D. Lind and  B. Marcus, {\sl An Introduction to Symbolic Dynamics and Coding,} 
Cambridge University Press, 1995. \ss

\item{[{\bf LS77}]} R.C. Lyndon and P.E. Schupp, Combinatorial Group Theory, Springer-Verlag,
Berlin, 1977.\ss

\item{[{\bf Na83}]} M. Nasu, Constant-to-one and onto global maps of homomorphisms between strongly connected graphs, {\sl Ergod.\ Th.\& Dynam.\ Sys. \bf 3} (1983), 387-411.\ss 

\item{[{\bf Qu95}]} F. Quinn, Lectures on axiomatic topological quantum field theory, {\sl Geometry and Quantum Field Theory} (D, Freed, K. Uhlenbeck, ed.), Amer.\ Math.\ Soc., 1995.\ss

\item{[{\bf Si93}]} D.S. Silver, Augmented group systems and $n$-knots, {\sl Mathematische
Annalen \bf 296} (1993), 585--593.\ss

\item{[{\bf SW96}]} D.S. Silver and S.G. Williams, Augmented group systems and shifts of finite
type, {\sl Israel\  J.\  Math.\ \bf 95} (1996), 231--251.\ss

\item{[{\bf SW99}]} D.S. Silver and S.G. Williams, Knot invariants from symbolic dynamical systems, {\sl Trans.\ Amer.\ Math.\ Soc.\ \bf 351} (1999), 3243--3265. \ss

\item{[{\bf SW99$'$}]} D.S. Silver and S.G. Williams, On groups with uncountably many subgroups of
finite index, {\sl J. Pure Appl. Algebra \bf 140} (1999), 75--86. \ss

\item{[{\bf SW04}]} D.S. Silver and S.G. Williams, Lifting representations of $\Z$-groups, preprint, April 2004.

\item{[{\bf St51}]} N. Steenrod, {\sl The Topology of Fibre Bundles,} Princeton Univ. Press, Princeton NJ, 1951. \ss

\item{[{\bf Wi73}]} R.F. Williams, Classification of subshifts of finite type, {\sl Annals
of Mathematics \bf 98} (1973), 120--153; erratum, {\sl Annals of Mathematics \bf 99} (1974),
380--381.\ss

\bs
\item{} Dept.\ of Mathematics and Statistics, Univ.\ of South Alabama, Mobile, AL  36688 
\item{} email: silver@jaguar1.usouthal.edu, swilliam@jaguar1.usouthal.edu

\end